\documentclass{amsart}
\usepackage[margin=1.0in]{geometry}
\usepackage{amsthm, amsmath, amssymb, mathtools, mathrsfs}
\usepackage{aliascnt}
\usepackage{hyperref}
\usepackage{breakurl}

\usepackage[textsize=footnotesize]{todonotes}

\usepackage{xcolor}
\definecolor{dblue}{rgb}{0,0,0.70}
\hypersetup{
	unicode=true,
	colorlinks=true,
	citecolor=dblue,
	linkcolor=dblue,
	anchorcolor=dblue
}

\makeatletter
\expandafter\g@addto@macro\csname th@plain\endcsname{%
	\thm@notefont{\bfseries}
}%
\expandafter\g@addto@macro\csname th@remark\endcsname{%
	\thm@headfont{\bfseries}
}%
\makeatother

\newtheorem{theorem}{Theorem}[section]	
\newtheorem*{theorem*}{Theorem}

\newaliascnt{lemma}{theorem}
\newtheorem{lemma}[lemma]{Lemma}
\aliascntresetthe{lemma}
\newtheorem*{lemma*}{Lemma}

\newaliascnt{proposition}{theorem}
\newtheorem{proposition}[proposition]{Proposition}
\aliascntresetthe{proposition}

\newaliascnt{corollary}{theorem}
\newtheorem{corollary}[corollary]{Corollary}
\aliascntresetthe{corollary}

\theoremstyle{remark}

\newaliascnt{remark}{theorem}

\aliascntresetthe{remark}
\newaliascnt{question}{theorem}
\newtheorem{question}[question]{Question}
\aliascntresetthe{question}

\newtheorem*{question*}{Question}

\newaliascnt{definition}{theorem}
\newtheorem{definition}[definition]{Definition}
\aliascntresetthe{definition}

\newaliascnt{example}{theorem}
\newtheorem{example}[example]{Example}
\aliascntresetthe{example}

\usepackage[backend=biber, bibstyle=numeric, citestyle=numeric]{biblatex}
\newcommand{\lr}{\leftrightarrow}
\newcommand{\Ord}{\mathrm{Ord}}

\newcommand{\dom}{\operatorname{dom}}

\newcommand{\rank}{\operatorname{rank}}
\newcommand{\crit}{\operatorname{crit}}
\newcommand{\restricts}{\mathrel{\upharpoonright}}

\newcommand{\TC}{\operatorname{TC}}

\newcommand{\mv}{\operatorname{mv}}
\newcommand{\mhyphen}{\text{-}}

\newcommand{\up}{\mathsf{up}}
\newcommand{\op}{\mathsf{op}}

\newcommand{\dnot}{{\lnot\lnot}}
\newcommand{\Low}{\operatorname{Low}}

\newcommand{\llbr}{[\mkern-2mu[}
\newcommand{\rrbr}{]\mkern-2mu]}
\newcommand{\lag}{\langle}
\newcommand{\rag}{\rangle}

\newcommand{\rrarrows}{\rightrightarrows}
\newcommand{\lrlrarrows}{\>\mathrlap{\leftleftarrows}{\rightrightarrows}\>}

\makeatletter
\def\oversortoftilde#1{\mathop{\vbox{\m@th\ialign{##\crcr\noalign{\kern3\p@}%
				\sortoftildefill\crcr\noalign{\kern3\p@\nointerlineskip}%
				$\hfil\displaystyle{#1}\hfil$\crcr}}}\limits}

\title{How strong is a Reinhardt set over extensions of CZF?}
\author{Hanul Jeon}
\email{ \href{mailto:hanuljeon95@gmail.com}{hanuljeon95@gmail.com}}
\urladdr{ \href{https://hanuljeon95.github.io}{https://hanuljeon95.github.io} }

\begin{document}
	\maketitle
	\begin{abstract}
	We investigate the lower bound of the consistency strength of $\mathsf{CZF}$ with Full Separation $\mathsf{Sep}$ and a Reinhardt set, a constructive analogue of Reinhardt cardinals.
	We show that $\mathsf{CZF+Sep}$ with a Reinhardt set interprets $\mathsf{ZF^-}$ with a cofinal elementary embedding $j:V\prec V$. We also see that $\mathsf{CZF+Sep}$ with a Reinhardt set interprets $\mathsf{ZF^-}$ with a model of $\mathsf{ZF+WA_0}$, the Wholeness axiom for bounded formulas.
\end{abstract}
	
	\section{Introduction}
	Large cardinals are one of the important topics of set theory: the linear hierarchy of large cardinals provides a scale to fathom the consistency strength of a given set-theoretic statement. Reinhardt introduced a quite strong notion of a large cardinal, now known as Reinhardt cardinal. Unfortunately, Kunen \cite{Kunen1971} showed that Reinhardt cardinals do not exist in $\mathsf{ZFC}$. However, Kunen's proof relied on the Axiom of Choice, so it remains the hope that Reinhardt cardinals are consistent if we dispose of the Axiom of Choice. Recently, there are attempts to study Reinhardt cardinals over $\mathsf{ZF}$ and find intrinsic evidence for \emph{choiceless large cardinals}, cardinals that are incompatible with the Axiom of Choice.
	
	We may ask choiceless large cardinals are actually consistent with subtheories of $\mathsf{ZFC}$, and there are some positive answers for this question. For example, Schultzenberg \cite{Schlutzenberg2020Kunen} showed that $\mathsf{ZF}$ with an elementary embedding $j:V_{\lambda+2}\prec V_{\lambda+2}$ is consistent if $\mathsf{ZFC}$ with $I_0$ is. Matthews \cite{Matthews2020} proved that Reinhardt cardinal is compatible with $\mathsf{ZFC^-}$ if we assume the consistency of $\mathsf{ZFC}$ with $I_1$.
	
	We may also ask the lower bound of the consistency strength of large cardinals over subtheories of $\mathsf{ZFC}$, which is quite non-trivial if we remove axioms other than choice. For example, Schultzenberg proved that $\mathsf{ZF}$ with $j:V_{\lambda+2}\prec V_{\lambda+2}$ is equiconsistent with $\mathsf{ZFC}$ with $I_0$.
	We do not have an actual bound for $\mathsf{ZFC^-}$ with $j:V\prec V$, but we can obtain the lower bound if we add the assumption $V_{\crit j}\in V$ that Matthews' model satisfies: in that case, we can see that $V_\lambda$ exists and it is a model of $\mathsf{ZFC}$ with the Wholeness axiom $\mathsf{WA}$.
	
	In this paper, we will take another `subtheory' of $\mathsf{ZFC}$ that lacks the law of excluded middle, namely $\mathsf{CZF}$.
	$\mathsf{CZF}$ is a weak theory in that its consistency strength is the same as that of Kripke-Platek set theory $\mathsf{KP}$ with Infinity. However, adding the law of excluded middle to $\mathsf{CZF}$ results in the full $\mathsf{ZF}$.
	The aim of this paper is to measure the `lower bound' of the consistency strength of $\mathsf{CZF}$ with Full Separation $\mathsf{Sep}$ and a Reinhardt set. Hence we have how hard to establish the consistency of $\mathsf{CZF+Sep}$ with a Reinhardt set. The main result of this paper is as follows, which is a consequence of \autoref{Theorem:CZFwithElementaryEmbeddings} and \autoref{Theorem:AnalysisEmbedding-main}.
	\begin{proposition}\label{Proposition:MainResult}
		The theory $\mathsf{CZF+Sep}$ + `there is a Reinhardt set' can interpret the following theory: $\mathsf{ZF^-}$ with the cofinal embedding $j:V\prec V$ with a transitive set $K$ such that $j(x)=x$ for all $x\in K$, $K\in j(K)$ and $\Lambda:=\bigcup_{n\in\omega}j^n(K)$ thinks it is a model of $\mathsf{ZF}$ with $\mathsf{WA}_0$, the Wholeness axiom for $\Delta_0$-formulas.
	\end{proposition}
	
	\subsection{The structure of this paper}
	In \autoref{Section:CZF} and \ref{Section:LargeSet}, we will cover relevant preliminaries. We will review constructive set theory in \autoref{Section:CZF} and large set axioms in \autoref{Section:LargeSet}, so the readers who are already familiar with these topics may skip them.
	In \autoref{Section:HeytingInterpretation}, we review and discuss Gambino's Heyting-valued interpretation \cite{Gambino2006} over the double negation formal topology $\Omega$, which is the main tool of the paper. It forces $\mathsf{\Delta_0\mhyphen LEM}$, the law of excluded middle for bounded formulas. By assuming the Full Separation, we may turn $\mathsf{\Delta_0\mhyphen LEM}$ to the full law of excluded middle.
	However, forcing over $\Omega$ does not preserve every axiom of $\mathsf{CZF}$: it does not preserve the Axiom of Subset Collection, a $\mathsf{CZF}$-analogue of Power Set axiom. Thus the resulting theory of forcing over $\Omega$ is $\mathsf{ZF^-}$ if we start from $\mathsf{CZF+Sep}$.
	Unlike `small' formal topologies, $\Omega$ is not absolute between transitive models, and it causes issues about the absoluteness of the Heyting-valued interpretation. We will also discuss it in this section.
	\autoref{Section:DNF-ElemEmbedding} is devoted to prove elementary embeddings are preserved under $\Omega$. However, $\Omega$ does not prove the inaccessibility of a critical point of the elementary embedding, and it restricts the analysis on the consistency strength of $\mathsf{CZF+Sep}$ with an elementary embedding. We will deal with this issue in \autoref{Section:CriticalPoint} by showing that the Heyting-valued interpretation under $\Omega$ still proves the critical point enjoys a strong reflection principle, which makes the critical point a transitive model of $\mathsf{ZF}$ with large cardinal axioms. We briefly discuss why achieving the upper bound of the consistency strength of our object theories in \autoref{Section:RemarkQuestions}, with some concluding questions.
	
	\section{Constructive set theory}\label{Section:CZF}
	In this section, we will briefly review $\mathsf{ZFC^-}$, $\mathsf{ZFC}$ without power set and the constructive set theory. There are various formulations of constructive set theories, but we will focus on $\mathsf{CZF}$.
	
	\subsection{$\mathsf{ZFC}$ without Power Set}
	We will frequently mention $\mathsf{ZFC}$ without Power Set, called $\mathsf{ZFC^-}$. However, $\mathsf{ZFC^-}$ is not obtained by just dropping Power Set from $\mathsf{ZFC}$: 
	\begin{definition}
		$\mathsf{ZFC^-}$ is the theory obtained from $\mathsf{ZFC}$ as follows: it drops Power Set, uses Collection instead of Replacement, and the well-ordering principle instead of the usual statement of Choice. $\mathsf{ZF^-}$ is a system obtained from $\mathsf{ZFC^-}$ by dropping the well-ordering principle. 
	\end{definition}
	Note that using Collection instead of Replacement is necessary to avoid pathologies. See \cite{GitmanHamkinsJohnstone2016} for the details. It is also known by \cite{FriedmanGitmanKanovei2019} that $\mathsf{ZFC^-}$ does not prove the reflection principle.
	
	\subsection{Axioms of $\mathsf{CZF}$}
	Constructive Zermelo-Fraenkel set theory $\mathsf{CZF}$ is introduced by Aczel \cite{Aczel1978} with his type-theoretic interpretation of $\mathsf{CZF}$. We will introduce subtheories called \emph{Basic Constructive Set Theory} $\mathsf{BCST}$ and $\mathsf{CZF^-}$ before defining the full $\mathsf{CZF}$.
	\begin{definition}
		$\mathsf{BCST}$ is the theory which consists of Extensionality, Pairnig, Union, Emptyset, Replacement, and $\Delta_0$-Separation.		
		$\mathsf{CZF^-}$ is obtained by adding the following axioms to $\mathsf{BCST}$: Infinity, $\in$-induction, and \emph{Strong Collection} which states the following: if $\phi(x,y)$ is a formula such that $\forall x\in a\exists y \phi(x,y)$ for given $a$, then we can find $b$ such that
		\begin{equation*}
			\forall x\in a\exists\in b y \phi(x,y) \land \forall y\in b\exists x\in a \phi(x,y).
		\end{equation*}
	\end{definition}
	
	We will also define some synonyms for frequently-mentioned axioms:
	\begin{definition}
		We will use $\mathsf{Sep}$, $\mathsf{\Delta_0\mhyphen Sep}$, $\mathsf{\Delta_0\mhyphen LEM}$ for denoting Full Separation (i.e., Separation for all formulas), $\Delta_0$-Separation and the Law of Excluded Middle for $\Delta_0$-formulas.
	\end{definition}
	
	Full separation proves Strong Collection from Collection, but $\Delta_0$-Separation is too weak to do it. It is also known that $\Delta_0$-Separation is equivalent to the existence of the intersection of two sets. See Section 9.5 of \cite{AczelRathjen2010} for its proof.
	\begin{proposition}\pushQED{\qed} 
		Working over $\mathsf{BCST}$ without $\Delta_0$-Separation, $\Delta_0$-Separation is equivalent to the \emph{Axiom of Binary Intersection}, which asserts that $a\cap b$ exists if $a$ and $b$ are sets. \qedhere
	\end{proposition}
	
	It is convenient to introduce the notion of \emph{multi-valued function} to describe Strong Collection and Subset Collection that appears later. Let $A$ and $B$ be classes. A relation $R\subseteq A\times B$ is a \emph{multi-valued function from $A$ to $B$} if $\dom R=A$. In this case, we write $R:A\rrarrows B$. We use the notation $R:A\lrlrarrows B$ if both $R:A\rrarrows B$ and $R:B\rrarrows A$ hold.
	Then we can rephrase Strong Collection as follows: for every set $a$ and a class $R:a\rrarrows V$, there is a set $b$ such that $R:a\rrarrows b$.
	
	Now we can state the Axiom of Subset Collection:
	\begin{definition}
		The Axiom of Subset Collection states the following: let $R_u$ be a class with a parameter $u\in V$. For each $a,b\in V$, we can find a set $c\in V$ such that
		\begin{equation*}
			R_u:a\rrarrows b\implies \exists d\in c (R_u:a\rrarrows d).
		\end{equation*}
		$\mathsf{CZF}$ is the theory by adding Subset Collection to $\mathsf{CZF^-}$.
	\end{definition}
	There is a simpler version of Subset Collection known as \emph{Fullness}, which is a bit easier to understand.
	\begin{definition}
		The Axiom of Fullness states the following: Let $\mv(a,b)$ the class of all multi-valued function from $a$ and $b$. Then there is a subset $c\subseteq \mv(a,b)$ such that if $r\in \mv(a,b)$, then there is $s\in c$ such that $s\subseteq r$.
		We call $c$ to be \emph{full} in $\mv(a,b)$.
	\end{definition}
	
	Then the following holds:
	\begin{proposition}\label{Proposition:SubsetCollection}
		\begin{enumerate}
			\item\label{Item:Fullness} \normalfont{($\mathsf{CZF^-}$)} Subset Collection is equivalent to Fullness.
			\item \normalfont{($\mathsf{CZF^-}$)} Power Set implies Subset Collection.
			\item \normalfont{($\mathsf{CZF^-}$)} Subset Collection proves the function set ${^a}b$ exists for all $a$ and $b$.
			\item \normalfont{($\mathsf{CZF^-}$)} If $\mathsf{\Delta_0\mhyphen LEM}$ holds, then Subset Collection implies Power Set.
		\end{enumerate}
	\end{proposition}
	We do not provide the proof for it, and the readers may consult with \cite{AczelRathjen2001} or \cite{AczelRathjen2010} for its proof. Note that Subset Collection does not increase the proof-theoretic strength of $\mathsf{CZF^-}$ while the Axiom of Power Set does.
	The following lemma is useful to establish \eqref{Item:Fullness} of \autoref{Proposition:SubsetCollection}, and is also useful to treat multi-valued functions:
	\begin{lemma}\label{Lemma:PrelimAdjuectmentFtn}
		Let $R: A\rightrightarrows B$ be a multi-valued function. Define 
		$\mathcal{A}(R) : A\rightrightarrows A\times B$ by
		\begin{equation*}
			\mathcal{A}(R) = \{\lag a,\lag a,b\rag\rag \mid 
			\lag a,b\rag \in R\},
		\end{equation*}
		then the following holds:
		\begin{enumerate}
			\item $\mathcal{A}(R) : A\rightrightarrows S\iff R\cap S:A\rightrightarrows B$,
			\item $\mathcal{A}(R) : A\leftleftarrows S\iff S\subseteq R$.
		\end{enumerate}
	\end{lemma}
	
	\begin{proof}
		For the first statement, observe that $\mathcal{A}(R): A\rightrightarrows S$ is equivalent to
		\begin{equation*}
			\forall a\in A \exists s\in S :\lag a,s\rag \in \mathcal{A}(R).
		\end{equation*}
		By the definition of $\mathcal{A}$, this is equivalent to
		\begin{equation*}\label{Formula:Mvaluedftn_eq00}
			\forall a\in A \exists s\in S [\exists b\in B: s=\lag a,b\rag \land \lag a,b\rag\in R].
		\end{equation*}
		We can see that the above statement is equivalent to $\forall a\in A\exists b\in B : \lag a,b\rag\in R\cap S$, which is the definition of $R\cap S : A\rightrightarrows B$.
		For the second claim, observe that $\mathcal{A}(R): A\leftleftarrows S$ is equivalent to
		\begin{equation*}
			\forall s\in S \exists a\in A :\lag a,s\rag \in \mathcal{A}(R).
		\end{equation*}
		By rewriting $\mathcal{A}$ to its definition, we have
		\begin{equation*}
			\forall s\in S \exists a\in A : [\exists b\in B : s=\lag a,b\rag \in R].
		\end{equation*}
		We can see that it is equivalent to $S\subseteq R$.
	\end{proof}
	
	The following lemma provides useful applications of Strong Collection:
	\begin{lemma}\label{Lemma:SetMV}
		If $a\in A$ and $R:a\rrarrows A$, then there is a set $b\in A$ such that $b\subseteq R$ and $b:a\rrarrows A$.
	\end{lemma}
	\begin{proof}
		Consider $\mathcal{A}(R):a\rrarrows a\times A$. By the second-order Strong Collection over $A$, there is $b\in A$ such that $\mathcal{A}(R):A\lrlrarrows b$. Hence by  \autoref{Lemma:PrelimAdjuectmentFtn}, we have $b\subseteq R$ and $b:a\rrarrows A$.
	\end{proof}
	
	\subsection{Inductive definition}
	Various recursive construction on $\mathsf{CZF}$ is given by inductive definition. The readers might refer \cite{AczelRathjen2001} or \cite{AczelRathjen2010} to see general information about inductive definition, but we will review some of it for the readers who are not familiar with it.
	\begin{definition}
		An \emph{inductive definition} $\Phi$ is a class of pairs $\lag X,a\rag$.
		For an inductive definition $\Phi$, associate $\Gamma_\Phi(C)=\{a\mid \exists X\subseteq C \lag X,a\rag\in\Phi\}$. A class $C$ is $\Phi$-closed if $\Gamma_\Phi(C)\subseteq C$.
	\end{definition}
	We may think $\Phi$ as a generalization of a deductive system, and $\Gamma_\Phi(C)$ a class of theorems derivable from the class of axioms $C$. Some authors use the notation $X\vdash_\Phi a$ or $X/a\in\Phi$ instead of $\lag X,a\rag\in \Phi$.
	
	Each Inductive definitions arise the least class fixed point:
	\begin{theorem}[Class Inductive Definition Theorem]
		Let $\Phi$ be an inductive definition. Then there is a smallest $\Phi$-closed class $I(\Phi)$.
	\end{theorem}
	
	The following lemma is the essential tool for the proof of Class Inductive Definition Theorem.
	See Lemma 12.1.2 of \cite{AczelRathjen2010} for its proof:
	\begin{lemma}\label{Lemma:ItreationClass}
		Every inductive definition $\Phi$ has a corresponding \emph{iteration class} $J$, which satisfies $J^a=\Gamma\left(\bigcup_{x\in a}J^x\right)$ for all $a$, where $J^a=\{x\mid \lag a,x\rag\in J\}$.
	\end{lemma}
	
	\subsection{$\mathsf{CZF}$ versus $\mathsf{IZF}$}
	There are two possible constructive formulations of $\mathsf{ZF}$, namely $\mathsf{IZF}$ and $\mathsf{CZF}$, although we will focus on the latter.
	\begin{definition}
		$\mathsf{IZF}$ is the theory that comprises the following axioms: Extensionality, Pairing, Union, Infinity, $\in$-induction, Separation, Collection, and Power Set.
	\end{definition}
	It is known that every theorem of $\mathsf{CZF}$ is that of $\mathsf{IZF}$. Moreover, $\mathsf{IZF}$ is quite strong in the sense that its proof-theoretic strength is the same as that of $\mathsf{ZF}$. On the other hand, it is known that the proof-theoretic strength of $\mathsf{CZF}$ is equal to that of Kripke-Platek set theory $\mathsf{KP}$ with Infinity.
	$\mathsf{IZF}$ is deemed to be \emph{impredicative} due to the presence of Full Separation and Power Set.\footnote{There is no consensus on the definition on the predicativity. The usual informal description of the predicativity is rejecting self-referencing definitions.} On the other hand, $\mathsf{CZF}$ is viewed as predicative since it allows the \emph{type-theoretic interpretation} given by Aczel \cite{Aczel1978}. However, adding the full law of excluded middle into $\mathsf{IZF}$ or $\mathsf{CZF}$ results in the same $\mathsf{ZF}$.
	
	\section{Large set axioms}\label{Section:LargeSet}
	In this section, we will discuss large set axioms, which is an analogue of large cardinal axioms over $\mathsf{CZF}$. Since ordinals over $\mathsf{CZF}$ could be badly behaved (for example, they need not be well-ordered), we focus on the structural properties of given sets to obtain higher infinities over $\mathsf{CZF}$.
	We also compare the relation between large cardinal axioms over well-known theories like $\mathsf{ZF}$ and large set axioms.
	
	\subsection{Tiny and Small Large set axioms}
	The first large set notions over $\mathsf{CZF}$ would be \emph{regular sets}. Regular sets appear first in Aczel's paper \cite{Aczel1986} about inductive definitions over $\mathsf{CZF}$. As we will see later, regular sets can `internalize' most inductive constructions, which turns out to be useful in many practical cases.
	
	\begin{definition}
		A transitive set $A$ is \emph{regular} if it satisfies second-order Strong Collection:
		\begin{equation*}
			\forall a\in A\forall R [R: a\rrarrows A\to \exists b\in A (R:a\lrlrarrows b)].
		\end{equation*}
		A regular set $A$ is \emph{$\bigcup$-regular} if $\bigcup a\in A$ for all $a\in A$.
		A regular set $A$ is \emph{inaccessible} if $(A,\in)$ is a model of $\mathsf{CZF}$, and furthermore, it also satisfies the second-order Subset Collection:
		\begin{equation*}
			\forall a,b\in A\exists c\in A \forall u\in A\forall R (R_u:a\rrarrows b)\to \exists d\in c (R_u:a\lrlrarrows d).
		\end{equation*}
		The \emph{Regular Extension Axiom} $\mathsf{REA}$ asserts that every set is contained in some regular set. The \emph{Inaccessible Extension Axiom} $\mathsf{IEA}$ asserts that every set is contained in an inaccessible set.
	\end{definition}
	There is no `pair-closed regular sets' since every regular set is closed under pairings if it contains 2:
	\begin{lemma}\label{Lemma:Regular2-Pairing}\pushQED{\qed}
		If $A$ is regular and $2\in A$, then $\lag a,b\rag \in A$ for all $a,b\in A$. \qedhere
	\end{lemma}
	$\mathsf{REA}$ has various consequences: For example, $\mathsf{CZF^-+REA}$ proves Subset Collection. Moreover, it also proves that every \emph{bounded} inductive definition $\Phi$ has a set-sized fixed point $I(\Phi)$.
	
	The notion of regular sets is quite restrictive, as it does not have Separation axioms, at least for $\Delta_0$-formulas, so we have no way to do any internal construction over a regular set.
	The following notion is a strengthening of regular set, which resolves the issue of internal construction:
	\begin{definition}
		A regular set $A$ is \emph{BCST-regular} if $A\models \mathsf{BCST}$. Equivalently, $A$ satisfies Union, Pairing, Empty set and Binary Intersection.
	\end{definition}
	We do not know that $\mathsf{CZF}$ proves every regular set is BCST-regular, although Lubarsky and Rathjen \cite{LubarskyRathjen2003} proved that the set of all hereditarily countable sets in the Feferman-Levy model is functionally regular but not $\bigcup$-regular. It is not even sure that the existence of a regular set implies that of BCST-regular set. However, every inaccessible set is BCST-regular, and every BCST-regular set appearing in this paper is inaccessible.

	What are regular sets and inaccessible sets in the classical context? The following result illustrates how these sets look like under the well-known classical context:
	\begin{proposition}
		\begin{enumerate}
			\item \normalfont{($\mathsf{ZF^-}$)} Every $\bigcup$-regular set containing $2$ is a transitive model of second-order $\mathsf{ZF^-}$, $\mathsf{ZF^-_2}$.
			
			\item \normalfont{($\mathsf{ZFC^-}$)} Every $\bigcup$-regular set containing $2$ is of the form $H_\kappa$ for some regular cardinal $\kappa$.
			
			\item \normalfont{($\mathsf{ZF^-}$)} Every inaccessible set is of the form $V_\kappa$ for some inaccessible $\kappa$.
		\end{enumerate}
	\end{proposition}
	Note that we follow the Hayut and Karagila's definition \cite{HayutKaragila2020} of inaccessiblity in choiceless context; that is, $\kappa$ is inaccessible if $V_\kappa\models \mathsf{ZF_2}$.
	\begin{proof}
		\begin{enumerate}
			\item Let $A$ be a regular set containing $2$. We know that $A$ satisfies Extensionality, $\in$-induction, Union and the second-order Collection. Hence it remains to show that the second-order Separation holds.
			
			Let $X\subseteq A$ and $a\in A$. Fix $c\in X\cap a$. Now consider the function $f:a\to A$ defined by
			\begin{equation*}
				f(x)=\begin{cases}
					x&\text{if }x\in X,\text{ and}\\
					c&\text{otherwise}.
				\end{cases}
			\end{equation*}
			By the second-order Strong Collection over $A$, we have $b\in A$ such that $f:a\lrlrarrows
			b$, and thus $b=a\cap X$.
			
			\item Let $A$ be a regular set. Let $\kappa$ be the least ordinal that is not a member of $A$. Then $\kappa$ must be a regular cardinal: if not, there is $\alpha<\kappa$ and a cofinal map $f:\alpha\to \kappa$. By transitivity of $A$ and the definition of $\kappa$, we have $\alpha\in A$, so $\kappa\in A$ by the second-order Replacement and Union, a contradiction.
			
			We can see that $\mathsf{ZFC^-}$ proves $H_\kappa$ is a class model of $\mathsf{ZFC^-}$, and $A$ satisfies the Well-ordering Principle. 
			We can also show that $A\subseteq H_\kappa=\{x:|\TC x|<\kappa\}$ holds: 
			We know that $A\cap\Ord= H_\kappa\cap \Ord=\kappa$. By the second-order Separation over $A$, $\mathcal{P}(\Ord)\cap A=\mathcal{P}(\Ord)\cap H_\kappa$. Hence $A=H_\kappa$: for each $x\in H_\kappa$, we can find $\theta<\kappa$, $R\subseteq \theta\times\theta$ and $X\subseteq\theta$ such that $(\operatorname{trcl} x,\in,x)\cong (\theta,R,X)$. (Here we treat $x$ as a unary relation.) Then $(\theta,R,X)\in A$, so $x\in A$ by Mostowski Collapsing Lemma.
			
			\item If $A$ is inaccessible, then $A$ is closed under the true power set of its elements, since the second-order Subset Collection implies if $a,b\in A$ then ${^a}b\in A$. Hence $A$ must be of the form $V_\kappa$ for some $\kappa$. Moreover, $\kappa$ is inaccessible because $V_\kappa=A\models \mathsf{ZF^-_2}$.
			\qedhere
		\end{enumerate}
	\end{proof}
	
	\begin{question}
		Is there a characterization of regular sets over $\mathsf{ZFC}$? How about $\bigcup$-regular sets over $\mathsf{ZF^-}$?
	\end{question}
	
	\subsection{Large Large set axioms}
	There is no reason to stop defining large set notions up to weaker ones. Hence we define stronger large set axioms. The main tool to access strong large cardinals (up to measurable cardinals) is to use elementary embedding, so we follow the same strategy:
	
	\begin{definition}[$\mathsf{CZF^-}$]
		Working over the extended language $\in$, a unary functional symbol $j$ and a unary predicate symbol $M$. We will extend $\mathsf{CZF^-}$ as follows: we allow $j$ and $M$ in $\Delta_0$-Separation and Strong Collection (also for Subset Collection if we start from $\mathsf{CZF}$), and add the following schemes:
		\begin{enumerate}
			\item $M$ is transitive, $\forall x M(j(x))$, and
			\item $\forall \vec{x} [\phi(\vec{x})\lr \phi^M(j(\vec{x}))]$ for every $\phi$ which does not contain any $j$ or $M$, where $\phi^M$ is the relativization of $\phi$ over $M$.
		\end{enumerate}
		If a transitive set $K$ satisfies $\forall x\in K j(x)=x$ and $K\in j(K)$, then we call $K$ a \emph{critical point} and we call $K$ a \emph{critical set} if $K$ is also inaccessible. If $M=V$ and $K$ is transitive and inaccessible, then we call $K$ a \emph{Reinhardt set}.
	\end{definition}
	We will use the term critical point and critical set simultaneously, so the readers should distinguish the difference of these two terms. (For example, a critical point need not be a critical set unless it is inaccessible.)
	Note that the definition of critical sets is apparently stronger than that is suggested by Hayut and Karagila \cite{HayutKaragila2020}. Finding the $\mathsf{CZF}$-definition of a critical set which is classically equivalent to a critical set in the style of Hayut and Karagila would be a good future work.
	Also, Ziegler \cite{ZieglerPhD} uses the term `measurable sets' to denote critical sets, but we will avoid this term for the following reasons: it does not reflect that the definition is given by an elementary embedding, and it could be confusing with measurable sets in measure theory.
	
	We do not know every elementary embedding $j:V\prec M$ over $\mathsf{CZF}$ enjoys cofinal properties. Surprisingly, the following lemma shows that $j$ become a cofinal map if $M=V$, even under $\mathsf{CZF}$ without any assumptions. 
	Note that the following lemma uses Subset Collection heavily. See Theorem 9.37 of \cite{ZieglerPhD} for its proof.
	\begin{lemma}[Ziegler \cite{ZieglerPhD}, $\mathsf{CZF}$]\label{Lemma:CZFReinhardtEmbeddingCofinal}
		Let $j:V\prec V$ be a non-trivial elementary embedding. Then $j$ is \emph{cofinal}, that is, we can find $y$ such that $x\in j(y)$ for each $x$. 
	\end{lemma}
	Note that Ziegler \cite{ZieglerPhD} uses the term \emph{set cofinality} to denote our notion of cofinality. However, we will use the term cofinality to harmonize the terminology with that of Matthews \cite{Matthews2020}.
	
	\section{Heyting-valued interpretation over the double negation formal topology}\label{Section:HeytingInterpretation}
	\subsection{General Heyting-valued interpretation}
	We will follow Gambino's definition \cite{Gambino2006} of Heyting-valued interpretation. We start this section by reviewing relevant facts on the Heyting-valued interpretation.
	
	Forcing is a powerful tool to construct a model of set theory. Gambino's definition of the Heyting-valued model (or alternatively, forcing) opens up a way to produce forcing models of $\mathsf{CZF^-}$. His Heyting-valued model starts from \emph{formal topology}, which formalizes a poset of open sets with a covering relation:
	\begin{definition}
		A structure $\mathcal{S}=(S,\le,\vartriangleleft)$ is \emph{formal topology} is a poset $(S,\le)$ endowed with $\vartriangleleft\subseteq S\times\mathcal{P}(S)$ such that
		\begin{enumerate}
			\item if $a\in p$, then $a\vartriangleleft p$,
			\item if $a\le b$ and $b\vartriangleleft p$, then $a\vartriangleleft p$,
			\item if $a\vartriangleleft p$ and $\forall x\in p (x\vartriangleleft q)$, then $a\vartriangleleft q$, and
			\item if $a\vartriangleleft p, q$, then $a\vartriangleleft (\downarrow p)\cap (\downarrow q)$, where $\downarrow p = \{r\in S\mid r\le p\}$.
		\end{enumerate}
	\end{definition}
	
	For each formal topology $\mathcal{S}$, we have a notion of \emph{nucleus} $\jmath p$ given by
	\begin{equation*}
		\jmath p = \{x\in S \mid x\vartriangleleft p\}.
	\end{equation*}
	Then the class $\Low(\mathcal{S})_\jmath$ of all lower subsets\footnote{A subset $p\subseteq S$ is a lower set if $\downarrow p = p$.} that are stable under $\jmath$ (i.e., $\jmath p = p$) form a \emph{set-generated frame}:
	\begin{definition}
		A structure $\mathcal{A}=(A,\le,\bigvee,\land,\top,g)$ is a \emph{set-generated frame} if $(A,\le,\bigvee,\land,\top)$ is a complete distributive lattice with the generating set $g\subseteq A$, such that the class $g_a=\{x\in g\mid x\le a\}$ is a set, and $a=\bigvee g_a$ for any $a\in A$.
	\end{definition}
	Note that every set-generated frame has every operation of Heyting algebra: for example, we can define $a\to b$ by $a\to b = \bigvee\{x\in g\mid x\land a\le b\}$, $\bot$ by $\varnothing$, and $\bigwedge p$ by $\bigvee\{x\in g\mid \forall y\in p (x\le y)\}$.
	
	\begin{proposition}
		For every formal topology $\mathcal{S}$, the class $\Low(\mathcal{S})_\jmath$ has a set-generated frame structure defined as follows:
		$p\land q=p\cap q$, $p\lor q=\jmath(p\cup q)$, $p\to q=\{x\in S\mid x\in p\to x\in q\}$, $\bigvee p=\jmath(\bigcup p)$, $\bigwedge p=\bigcap p$, $\top = S$, $\le$ as the inclusion relation, and $g=\{\{x\}\mid x\in S\}$.
	\end{proposition}
	We extend the nucleus $\jmath$ to general classes by taking $JP:=\bigcup\{\jmath p \mid p\subseteq P\}$, and define operations on classes by $P\land Q=P\cap Q$, $P\lor Q=J(P\cup Q)$ and $P\to Q =\{x\in S\mid x\in P\to x\in Q\}$. For a set-indexed collection of classes $\{P_x\mid x\in I\}$, take $\bigwedge_{x\in I} P_x=\bigcap_{x\in I} P_x$ and $\bigvee_{x\in I} P_x=J\left(\bigcup_{x\in I} P_x\right)$. Note that $JP=\jmath P$ if $P$ is a set and we have the Full Separation.
	
	The Heyting universe $V^\mathcal{S}$ over $\mathcal{S}$ is defined inductively as follows: $a\in V^\mathcal{S}$ if and only if $a$ is a function from $\dom a\subseteq V^\mathcal{S}$ to $\Low(\mathcal{S})_\jmath$. For each set $x$, we have the canonical representation $\check{x}$ of $x$ defined by $\dom \check{x}=\{\check{y}\mid y\in x\}$ and $\check{x}(\check{y})=\top$.
	Then define the Heyting interpretation $\llbr\phi\rrbr$ with parameters of $V^\mathcal{S}$ as follows:
	\begin{itemize}
		\item $\llbr a= b\rrbr=\left(\bigwedge_{x\in\dom a}a(x)\to\bigvee_{y\in\dom b} b(y)\land\llbr x=y\rrbr\right) \land \left(\bigwedge_{y\in\dom b}b(y)\to\bigvee_{x\in\dom a}a(x)\land\llbr x=y\rrbr\right)$,
		\item $\llbr a\in b\rrbr=\bigvee_{y\in\dom b} b(y)\land\llbr a=y\rrbr$,
		\item $\llbr\bot\rrbr=\bot$, $\llbr \phi\land \psi\rrbr=\llbr\phi\rrbr\land \llbr\psi\rrbr$, $\llbr \phi\lor \psi\rrbr=\llbr\phi\rrbr\lor \llbr\psi\rrbr$, and $\llbr \phi\to \psi\rrbr=\llbr\phi\rrbr\to \llbr\psi\rrbr$,
		\item $\llbr\forall x\in a\phi(x)\rrbr=\bigwedge_{x\in\dom a}a(x)\to \llbr\phi(x)\rrbr$ and $\llbr\exists x\in a\phi(x)\rrbr=\bigvee_{x\in\dom a}a(x)\land \llbr\phi(x)\rrbr$,
		\item $\llbr\forall x\phi(x)\rrbr=\bigwedge_{x\in V^\mathcal{S}}\llbr\phi(x)\rrbr$ and $\llbr\exists x\phi(x)\rrbr=\bigvee_{x\in V^\mathcal{S}}\llbr\phi(x)\rrbr$.
	\end{itemize}
	Then the interpretation validates every axiom of $\mathsf{CZF^-}$ and more:
	\begin{theorem}\label{Theorem:CZFPersistence}
		Working over $\mathsf{CZF^-}$, the Heyting-valued model $V^\mathcal{S}$ also satisfies $\mathsf{CZF^-}$. If $\mathcal{S}$ is set-presented and Subset Collection holds, then $V^\mathcal{S}\models \mathsf{CZF}$.
		If our background theory satisfies Full Separation, then so does $V^\mathcal{S}$.
	\end{theorem}
	\begin{proof}
		The first part of the theorem is shown by \cite{Gambino2006}, so we will concentrate on the preservation of Full Separation.
		For Full Separation, it suffices to see that the proof for bounded separation over $V^\mathcal{S}$ also works for Full Separation, since Full Separation ensures $\llbr\theta\rrbr$ is a set for every formula $\theta$.
	\end{proof}
	Let us finish this subsection with some constructors, which we need in a later proof.
	\begin{definition}
		For $\mathcal{S}$-names $a$ and $b$, $\up(a,b)$ is defined by $\dom(\up(a,b))=\{a,b\}$ and $(\up(a,b))(x)=\top$. $\op(a,b)$ is the name defined by $\op(a,b)=\up(\up(a,a),\up(a,b))$
	\end{definition}
	$\up(a,b)$ represents the unordeded pair $\{a,b\}$ over $V^\mathcal{S}$. Hence the name $\op(a,b)$ represents the ordered pair given by $a$ and $b$ over $V^\mathcal{S}$.
	
	\subsection{Double negation formal topology}
	Our main tool in this paper is the Heyting-valued interpretation with the \emph{double negation formal topology}. Unlike set-sized realizability or set-represented formal topology, the double negation topology and the resulting Heyting-valued interpretation need not be absolute between BCST-regular sets or transitive models of $\mathsf{CZF^-}$. Hence we need a careful analysis of the double negation formal topology, which is the aim of this subsection.
	
	\begin{definition}
		The \emph{double negation formal topology} $\Omega$ is the formal topology $(1,=,\vartriangleleft)$, where $x\vartriangleleft p$ if and only if $\lnot\lnot(x\in p)$.
	\end{definition}
	We can see that the class of lower sets $\Low(\Omega)$ is just the power set of 1, and the nucleus of $\mathcal{S}$ is given by the double complement 
	\begin{equation*}
		p^\dnot=\{0\mid \lnot\lnot(0\in p)\}.
	\end{equation*}
	Hence the elements of $\Low(\Omega)_\jmath$ is the collection of all \emph{stable} subsets of 1, that is, a set $p\subseteq 1$ such that $p=p^\dnot$.

	We will frequently mention the relativized Heyting-valued interpretation, the definition of relativized one is not different from the usual $V^\Omega$ and $\llbr\cdot\rrbr$, thus we do not introduce its definition. Notwithstanding that, it is still worth mentioning the notational convention for relativization:
	\begin{definition}
		Let $A$ be a transitive model of $\mathsf{CZF^-}$. Then $A^\Omega:=(V^\Omega)^A$ is the relativized $\Omega$-valued universe to $A$. If $A$ is a set, then $\tilde{A}$ denotes the $\Omega$-name defined by $\dom \tilde{A}:=A^\Omega$ and $\tilde{A}(x)=\top$ for all $x\in\dom \tilde{A}$.
	\end{definition}
	Note that the definition of $\tilde{A}$ makes sense due to \autoref{Lemma:HeytingUniversePreserving}. We often confuse $\tilde{A}$ and $A^\Omega$ if context is clear. It is also worth to mention that if $j$ is an elementary embedding, then $j(\tilde{K})=\widetilde{j(K)}$, so we may write $j^n(\tilde{K})$ instead of $\widetilde{j^n(K)}$.
	
	We cannot expect that $\Low(\Omega)$ is absolute between transitive models of $\mathsf{CZF^-}$, and as a result, we do not know whether its Heyting-valued universe $V^\Omega$ and Heyting-valued interpretation $\llbr\cdot\rrbr$ is absolute. Fortunately, the formula $p\in\Low(\Omega)_\jmath$ is $\Delta_0$, so it is absolute between transitive models of $\mathsf{CZF^-}$. As a result, we have the following absoluteness result on the Heyting-valued universe:	
	
	\begin{lemma}\label{Lemma:HeytingUniversePreserving}
		Let $A$ be a transitive model of $\mathsf{CZF^-}$ without Infinity. Then we have $A^\Omega=V^\Omega\cap A$. Moreover, if $A$ is a set, then $A^\Omega$ is also a set.
	\end{lemma}
	\begin{proof}
		We will follow the proof of Lemma 6.1 of \cite{Rathjen2003Realizability}.
		Let $\Phi$ be the inductive definition given by
		\begin{equation*}
			\lag X,a \rag\in\Phi\iff \text{$a$ is a function such that $\dom a\subseteq X$ and $a(x)\subseteq 1$, $a(x)^\dnot=a(x)$ for all $x\in\dom a$}.
		\end{equation*}
		We can see that the $\Phi$ defines the class $V^\Omega$. Furthermore, $\Phi$ is $\Delta_0$, so it is absolute between transitive models of $\mathsf{CZF^-}$.
		By Lemma \ref{Lemma:ItreationClass}, we have a class $J$ such that $V^\Omega = \bigcup_{a\in V} J^a$, and for each $s\in V$, $J^s=\Gamma_\Phi(\bigcup_{t\in s} J^t)$.
		Now consider the operation $\Upsilon$ given by
		\begin{equation*}
			\Upsilon(X):=\{a\in A\mid\exists Y\in A(Y\subseteq X\land \lag Y,a\rag\in\Phi)\}.
		\end{equation*}
		By Lemma \ref{Lemma:ItreationClass} again, there is a class $Y$ such that $Y^s=\Upsilon(\bigcup_{t\in s}Y^t)$ for all $s\in V$.
		Furthermore, we can see that $Y^s\subseteq V^\Omega$ by induction on $a$.
		
		Let $Y=\bigcup_{s\in A}Y^s$. We claim by induction on $s$ that $J^s\cap A\subseteq Y$.
		Assume that $J^t\cap A\subseteq Y$ holds for all $t\in s$. If $a\in J^s\cap A$, then the domain of $a$ is a subset of $A\cap \left(\bigcup_{t\in s} J^s\right)$, which is a subclass of $Y$ by the inductive assumption and the transitivity of $A$.
		Moreover, for each $x\in\dom a$ there is $u$ such that $x\in Y^u$.
		By Strong Collection over $A$, there is $v\in A$ such that for each $x\in \dom a$ there is $u\in v$ such that $x\in Y^u$. Hence $\dom a\subseteq \bigcup_{u\in v} Y^u$, which implies $a\in Y^v\subseteq Y$. 
		
		Hence $V^\Omega\cap A\subseteq Y$, and we have $Y=V^\Omega\cap A$. We can see that the construction of $Y$ is the relativized construction of $V^\Omega$ to $A$, so $Y=A^\Omega$. Hence $A^\Omega=V^\Omega\cap A$.
		If $A$ is a set, then $\Upsilon(X)$ is a set for each set $X$, so we can see by induction on $a$ that $Y^a$ is also a set for each $a\in A$. Hence $A^\Omega=Y=\bigcup_{a\in A} Y^a$ is also a set.
	\end{proof}
	
	We extended nucleus $\jmath$ to $J$ for subclasses of $\Low\mathcal{S}$, and use it to define the validity of formulas of the forcing language.
	We are working with the specific formal topology $\mathcal{S}=\Omega$, and in that case, $JP$ for a class $P\subseteq 1$ coincides with $JP = \bigcup\{q^\dnot \mid q\subseteq P\}$.
	It is easy to see that $P\subseteq JP\subseteq P^\dnot$. We also define the following relativized notion for any transitive class $A$ such that $1\in A$:
	\begin{equation*}
		J^AP = \bigcup\{q^\dnot \mid q\subseteq P\text{ and }q\in A\}.
	\end{equation*}
	If $P\in A$, then $J^AP=P^\dnot$, and in general, we have $P\subseteq J^AP\subseteq JP\subseteq P^\dnot$. Moreover, we have
	\begin{lemma}\label{Lemma:PropertiesOfRelativizedJ}
		Let $A$ and $B$ be transitive classes such that $1\in A,B$ and $P\subseteq 1$ be a class.
		\begin{enumerate}
			\item $A\subseteq B$ implies $J^AP\subseteq J^BP$.
			\item If $\mathcal{P}(1)\cap A=\mathcal{P}(1)\cap B$, then $J^AP=J^BP$.
		\end{enumerate}
	\end{lemma}
	
	However, the following proposition shows that we cannot prove they are the same from $\mathsf{CZF^-}$:
	\begin{proposition}
		\begin{enumerate}
			\item If $A\cap \mathcal{P}(1)=2$ (it holds when $A=2$ or $A=V$ and $\mathsf{\Delta_0\mhyphen LEM}$ holds), then $J^AP=P$.
			
			\item If $P^\dnot\subseteq JP$ for every class $P$, then $\mathsf{\Delta_0\mhyphen LEM}$ implies the law of excluded middle for arbitrary formulas. \qed
		\end{enumerate}
	\end{proposition}
	
	$J^A$ has a crucial role in defining Heyting-valued interpretation, but it could differ from transitive set to transitive set. It causes absoluteness problems, which is apparently impossible to emend in general. The following lemma states facts on relativized Heyting interpretations:
	
	\begin{lemma}\label{Lemma:HeytingInterpretationBetweenModels}
		Let $A\subseteq B$ be transitive models of $\mathsf{CZF^-}$. Assume that $\phi$ is a formula with parameters in $A^\Omega$.
		\begin{enumerate}
			\item If $\phi$ is bounded, then $\llbr\phi\rrbr^A=\llbr\phi\rrbr^B$.
			\item If $\phi$ only contains bounded quantifications, logical connectives between bounded formulas, unbounded $\forall$, and $\land$, then $\llbr\phi\rrbr^A=\llbr\phi^{\tilde{A}}\rrbr^B$.
			\item If every conditional of $\to$ appearing in $\phi$ is bounded, then $\llbr\phi\rrbr^A\subseteq \llbr\phi^{\tilde{A}}\rrbr^B$.
			\item If $\mathcal{P}(1)\cap A=\mathcal{P}(1)\cap B$, then $\llbr\phi\rrbr^A=\llbr\phi^{\tilde{A}}\rrbr^B$.
		\end{enumerate}
	\end{lemma}
	\begin{proof}
		If $\phi$ is bounded, then $\llbr\phi\rrbr$ is defined in terms of double complement, Heyting connectives between subsets of 1, and set-sized union and intersection. These notions are absolute between transitive sets, so we can prove $\llbr\phi\rrbr$ is also absolute by induction on $\phi$. (In the case of atomic formulas, we apply the induction on $A^\Omega$-names.)
		
		The remaining clauses follow from the induction on $\phi$:
		For the unbounded $\forall$, we have
		\begin{equation*}
			\llbr\forall x\phi(x)\rrbr^A = \bigwedge_{x\in A^\Omega}\llbr \phi(x)\rrbr^A
			\subseteq \bigwedge_{x\in A^\Omega} \llbr\phi^{\tilde{A}}(x)\rrbr^B
			= \llbr\forall x\in\tilde{A} \phi^{\tilde{A}}(x)\rrbr^B.
		\end{equation*}
		under conditions in the remaining clauses.
		The case for $\land$ and $\to$ are similar. In the case of the second or fourth clause, we can see that the above argument raises equality.
		
		For the unbounded $\exists$, we have
		\begin{equation*}
			\llbr\exists x\phi(x)\rrbr^A = J^A \left(\bigcup  \{\llbr \phi(x)\rrbr^A\mid x\in A^\Omega\}\right)
			\subseteq J^B \left(\bigcup  \{\llbr \phi^{\tilde{A}}(x)\rrbr^B\mid x\in A^\Omega\}\right) = \llbr\exists x\in \tilde{A} \phi^{\tilde{A}}(x) \rrbr^B.
		\end{equation*}
		Note that if $A\in B$, then $\bigcup_{x\in A^\Omega} \llbr \phi^{\tilde{A}}(x)\rrbr^B\in \mathcal{P}(1)\cap B$, 
		thus $J^B\left(\bigcup_{x\in A^\Omega} \llbr \phi^{\tilde{A}}(x)\rrbr^B\right) = \left(\bigcup_{x\in A^\Omega} \llbr \phi^{\tilde{A}}(x)\rrbr^B\right)^\dnot$ in this case.
		
		The case for $\lor$ is similar to that of the unbounded $\exists$. For the last clause, we need $J^A=J^B$ that follows from $\mathcal{P}(1)\cap A=\mathcal{P}(1)\cap B$.
	\end{proof}
	
	As a final remark, note that we may understand the forcing over $\Omega$ as a double negation translation \`a la forcing. See 2.3 of \cite{Grayson1979} or \cite{Bell2014} for details.
	
	\section{Double negation translation of an elementary embedding}\label{Section:DNF-ElemEmbedding}
	The following section is devoted to the following theorem:
	\begin{theorem}\label{Theorem:MainOfDNFcritical} $(\mathsf{CZF+Sep})$
		Let $j:V\prec V$ be an elementary embedding and let $K$ be an inaccessible critical point of $j$; that is, $K$ is inaccessible, $K\in j(K)$ and $j(x)=x$ for all $x\in K$. Then there is a Heyting-valued model $V^\Omega$ such that $V^\Omega\models\mathsf{ZF^-}$ and it thinks $j:V\prec V$ is cofinal and has a critical point.
	\end{theorem}
	
	We will mostly follow the proof of Ziegler \cite{ZieglerPhD}, but we need to check his proof works on our setting since his applicative topology does not cover Heyting algebra generated by formal topologies that are not set-presentable.
	We will focus on Reinhardt sets, so considering the target universe $M$ of $j$ might be unnecessary. Nevertheless, we will consider $M$ and we will not assume $M=V$ unless it is necessary.
	
	We need to redefine Heyting-valued interpretations to handle critical sets and Reinhardt sets, so we define the interpretation of $M$ and $j$ in the forcing language.
	Since $j$ preserves names, we can interpret $j$ as $j$ itself. We will interpret $M$ as the $M^\Omega$. Thus, for example, $\llbr\forall x\in M \phi(x)\rrbr=\bigwedge_{x\in M^\Omega} \llbr\phi(x)\rrbr$.
	We discussed that the Heyting interpretation $\llbr\phi\rrbr$ need not be absolute between transitive sets. A cacophony of absoluteness issues causes technical trouble in an actual proof, so we want to assure $\llbr\phi\rrbr=\llbr\phi\rrbr^M$, which follows from $\mathcal{P}(1)=\mathcal{P}(1)\cap M$.
	The following lemma is due to Ziegler, and see Remark 9.51 to 9.52 of \cite{ZieglerPhD} for its proof:
	\begin{lemma}\label{Lemma:PowerSetPreserving}
		Let $p\in\mathcal{P}(1)$ 
		and $j:V\prec M$ be an elementary embedding. Then $j(p)=p$. Especially, we have $\mathcal{P}(1)=\mathcal{P}(1)\cap M$.
	\end{lemma}
	Hence by \autoref{Lemma:HeytingInterpretationBetweenModels}, $\llbr \phi\rrbr^M = \llbr \phi^{M^\Omega}\rrbr$ for any formula $\phi$ with parameters in $M^\Omega$.
	Thus we do not need to worry about the absoluteness issue on the Heyting interpretation.
	
	We are ready to extend our forcing language to $\{\in, j, M\}$. For each formula $\phi$, define
	\begin{equation*}
		\llbr\forall x\in M \phi(x)\rrbr := \bigwedge_{x\in M^\Omega} \llbr \phi(x)\rrbr
		\qquad\text{and}\qquad
		\llbr\exists x\in M \phi(x)\rrbr := \bigvee_{x\in M^\Omega} \llbr \phi(x)\rrbr
	\end{equation*}
	and define the remaining as given by \cite{Gambino2006}. (We may define $x\in M$ by using that this is equivalent to $\exists y\in M (x=y)$.)
	From this definition, we have an analogue of Lemma 4.26 of \cite{ZieglerPhD}, which is useful to check that $j$ is still elementary over $V^\Omega$:
	\begin{lemma}\label{Lemma:ValuePreserving}
		For any bounded formula $\phi(\vec{x})$ with all free variables displayed in the language $\in$ (that is, without $j$ and $M$), we have
		\begin{equation*}
			\llbr \phi(\vec{a}) \rrbr = \llbr \phi^{M^\Omega} (j(\vec{a})) \rrbr = \llbr \phi (j(\vec{a})) \rrbr
		\end{equation*}
		for every $\vec{a}\in V^\Omega$.
	\end{lemma}
	\begin{proof}
		For the first equality, we have
		\begin{equation}\label{Formula:Reinhardt000}
			\llbr \phi(\vec{a}) \rrbr = j(\llbr \phi(\vec{a}) \rrbr)
			= \llbr \phi^{M^\Omega}(j(\vec{a})) \rrbr.
		\end{equation}
		by \autoref{Lemma:PowerSetPreserving}. Note that the last equality of \eqref{Formula:Reinhardt000} follows from the induction on $\phi$. The second equality also follows from the induction on $\phi$.
	\end{proof}
	Moreover, we can check the following equalities easily:
	\begin{proposition}
		\begin{enumerate}
			\item $\llbr \forall x,y (x=y\to j(x)=j(y))\rrbr=\top$,
			\item $\llbr \forall x j(x)\in M\rrbr=\top$,
			\item $\llbr \forall x (x\in M\to \forall y\in x(y\in M)) \rrbr=\top$.
		\end{enumerate}
	\end{proposition}
	\begin{proof}
		The first equality follows from $\llbr x=y\rrbr = \llbr j(x)=j(y)\rrbr$, and the remaining two follow from the direct calculation.
	\end{proof}
	
	\begin{lemma}
		For every $\vec{a}\in V^\Omega$ and a formula $\phi$ that does not contain $j$ or $M$, we have $\llbr \phi(\vec{a})\lr \phi^M(j(\vec{a})) \rrbr=\top$.
	\end{lemma}
	\begin{proof}
		\autoref{Lemma:ValuePreserving} shows that this lemma holds for bounded formulas $\phi$.
		We will use full induction on $\phi$ to prove $\llbr\phi(\vec{a})\rrbr=\llbr\phi^M(j(\vec{a}))\rrbr$ for all $\vec{a}\in V^\Omega$.
		If $\phi$ is $\forall x\psi(x,\vec{a})$, we have 
		$\llbr \forall x\phi(x,\vec{a})\rrbr=\bigwedge_{x\in V^\Omega}\llbr\phi(x,\vec{a})\rrbr$.
		Since
		\begin{align*}
			0\in \bigwedge_{x\in V^\Omega}\llbr \phi(x,\vec{a})\rrbr &\iff \forall x\in V^\Omega (0\in \llbr\phi(x,\vec{a})\rrbr)\\
			&\iff \forall x \in (V^\Omega)^M (0\in \llbr\phi(x,j(\vec{a}))\rrbr),
		\end{align*}
		where the last equivalence follows from applying $j$ to the above formula. Since the last formula is equivalent to $0\in \bigwedge_{x\in M^\Omega} \llbr\phi(x)\rrbr$, we have 
		$\bigwedge_{x\in V^\Omega} \llbr\phi(x)\rrbr=\bigwedge_{x\in M^\Omega} \llbr\phi(x)\rrbr=\llbr \forall x\in M \phi^M(x)\rrbr$.
		
		If $\phi$ is $\exists x\psi(x,a)$, we have 
		$\llbr \exists x\phi(x,a)\rrbr=\bigvee_{x\in V^\Omega}\llbr\phi(x,\vec{a})\rrbr$. Moreover,
		\begin{align*}
			0\in \bigvee_{x\in V^\Omega}\llbr \phi(x,\vec{a})\rrbr 
			&\iff \exists p\subseteq 1 \left[ p\subseteq \bigcup_{x\in V^\Omega} \llbr\phi(x,\vec{a})\rrbr \text{ and } 0\in p^\dnot \right]
			\\&\iff \exists p\subseteq 1 \left[ p\subseteq \bigcup_{x\in (V^\Omega)^M} \llbr\phi^M(x,j(\vec{a}))\rrbr \text{ and } 0\in p^\dnot \right]
		\end{align*}
		Hence $0\in \bigvee_{x\in V^\Omega}\llbr \phi(x,\vec{a})\rrbr$ if and only if $0\in \bigvee_{x\in M^\Omega}\llbr \phi(x,j(\vec{a}))\rrbr$
	\end{proof}
	
	We extended the language of set theory to treat elementary embedding $j$ and the target universe $M$. We also expand the Full Separation $\mathsf{Sep}$ and Strong Collection to the extended language. Thus we need to check that Full Separation and Strong Collection under the extended language are also persistent under the double negation interpretation.
	We can see that the proof given by \cite{Gambino2006} and \autoref{Theorem:CZFPersistence} carries over, so we have the following:
	\begin{proposition}
		\begin{enumerate}
			\item If we assume Full Separation for the extended language, then $V^\Omega$ also thinks Full Separation for the extended language holds.
			\item $V^\Omega$ thinks Strong Collection for the extended language holds.
		\end{enumerate}
	\end{proposition}
	
	The essential property of a critical set is that it is inaccessible. Unfortunately, there is no hope to preserve the inaccessibility of a critical set. The main reason is that an inaccessible set must satisfy second-order Set Collection, which is not preserved by forcing under $\Omega$.
	Fortunately, being a critical point is preserved provided if it is regular:
	\begin{lemma}
		Let $K$ be a regular set such that $K\in j(K)$ and $j(x)=x$ for all $x\in K$.
		Then $\llbr\check{K}\in j(\tilde{K})\land \forall x\in \tilde{K}(j(x)=x)\rrbr=\top$.
	\end{lemma}
	\begin{proof}
		Since $(\text{$j(K)$ is inaccessible})^M$, $j(K)^\Omega = (j(K)^\Omega)^M = j(K)\cap V^\Omega$ is a set by \autoref{Lemma:HeytingUniversePreserving}.
		Also, $K\in j(K)$ implies $\tilde{K}\in j(K)$. Since the domain of $j(\tilde{K})=\widetilde{j(K)}$ is $j(K)\cap V^\Omega$, we have $\tilde{K}\in \dom j(\tilde{K})$, which implies $\llbr\tilde{K}\in j(\tilde{K})\rrbr=\top$.
		For the second assertion, observe that if $x\in \dom\tilde{K}$ then $j(x)=x$, so we have the desired conclusion.
	\end{proof}
	Note that the canonical name $\check{K}$, is also a critical point in $V^\Omega$. However, we stick to use $\tilde{K}$, whose reason becomes apparent in \autoref{Section:CriticalPoint}.
	
	We do not make use of Full Separation or Subset Collection until now. However, the following proof requires Full Separation (or at least, Separation for $\Sigma$-formulas) or $\mathsf{REA}$.
	\begin{lemma}
		If $j:V\prec V$ is a cofinal elementary embedding, then $V^\Omega$ thinks $j$ is cofinal.
	\end{lemma}
	\begin{proof}
		Let $a\in V^\Omega$. Then there is a set $X$ such that $a\in j(X)$. If we assume Full Separation, then $X\cap V^\Omega$ is a set and $j(X\cap V^\Omega)=j(X)\cap V^\Omega$.
		Let $b$ be a name such that $\dom b= X\cap V^\Omega$ and $b(y)=1$ for all $y\in\dom b$.
		Then $\llbr a\in j(b) \rrbr=\top$.
		
		We need some work if we assume $\mathsf{REA}$ instead: Take a set $X$ such that $a\in j(X)$. By $\mathsf{REA}$, we can find a regular $Y$ such that $X\in Y$. By \autoref{Lemma:HeytingUniversePreserving}, $Y^\Omega=Y\cap V^\Omega$ is a set. The remaining argument is identical to the previous one.
	\end{proof}
	
	Combining the above lemmas shows our main theorem and more:
	\begin{theorem}\label{Theorem:CZFwithElementaryEmbeddings}
		The consistency of the former implies the consistency of the latter. Here the former always assume that the critical point of $j$ is regular (but not for the latter), $M$ can be equal to $V$, and we always assume that $j$ is non-trivial:
		\begin{enumerate}
			\item $\mathsf{CZF^-}+j:V\prec M$ and $\mathsf{CZF^-+\Delta_0\mhyphen LEM}+j:V\prec M$,
			\item $\mathsf{CZF+REA}+j:V\prec V$ and $\mathsf{CZF^-+\Delta_0\mhyphen LEM}+j:V\prec V$ is cofinal,
			\item $\mathsf{CZF^-+Sep} + j:V\prec V\text{ (is cofinal)}$ and $\mathsf{ZF^-} + j:V\prec V\text{ (is cofinal)}$,
			\item $\mathsf{CZF+Sep}+j:V\prec V$ and $\mathsf{ZF^-}+j:V\prec V$ is cofinal,
			\item $\mathsf{IZF}+j:V\prec M$ and $\mathsf{ZF}+j:V\prec M$,
		\end{enumerate}
	\end{theorem}
	\begin{proof}
		It follows from our previous lemmas in this section, \autoref{Lemma:CZFReinhardtEmbeddingCofinal}, and that Heyting-valued interpretation preserves $\mathsf{Sep}$ and $\mathsf{Pow}$.
	\end{proof}
	
	Unfortunately, the above result answers little about the consistency strength of $\mathsf{CZF+Sep}$ with a Reinhardt set, since we know little about the consistency strength of $\mathsf{ZF^-}$ with a cofinal embedding $j:V\prec V$.
	The author thought that its consistency strength is very high from the following argument: consider the principle known as the \emph{Relation Reflection Scheme} ($\mathsf{RRS}$) defined by Aczel \cite{Aczel2008}. It is known that $\mathsf{RRS}$ is persistent under Heyting-valued interpretation, and $\mathsf{ZFC^-}$ proves the equivalence between $\mathsf{RRS}$ and the reflection principle. 
	Moreover, in an earlier version of \cite{Matthews2020}, Matthews claimed that $\mathsf{ZFC^-}$ with the reflection principle and the existence of a cofinal elementary embedding is sufficient to derive the contradiction. Hence Reinhardt cardinals are incompatible with choice over $\mathsf{ZF^-}$ with the reflection principle.
	However, Matthews pointed out to the author that the claim as stated above is incorrect: the problem is that the reflection principle over $\mathsf{ZFC^-}$ does not ensure that $V_{\crit j}$ is a set, which seems necessary to derive the contradiction.
	In spite of that, we can still see that assuming $\mathsf{DC}_\mu$-scheme for every cardinal $\mu$ proves there is no cofinal elementary embedding over $\mathsf{ZFC^-}$.
	
	Can we see the revised result as a form of Kunen's inconsistency phenomenon, so that it is evidence of the consistency strength of $\mathsf{ZF^-}$ with a cofinal elementary embedding? 
	We can still prove that $\mathsf{ZFC^-}$ and $\mathsf{ZFC^-}$ with $\mathsf{DC}_\mu$-scheme for all $\mu$ are equiconsistent since $L$ satisfies global choice. We know that $L$ is compatible with small large cardinals over $\mathsf{ZFC}$, so we may guess the same holds for $\mathsf{ZFC^-}$. However, we do not aware well about large cardinals over $\mathsf{ZFC^-}$ to conclude that the consistency strength of $\mathsf{ZF^-}$ with a cofinal elementary is high. Even worse, $L$ is not compatible with large large cardinals. Thus we cannot extend the same argument further.
	Analyzing the notion of large cardinals over $\mathsf{ZFC^-}$ and their consistency with $\mathsf{DC}_\mu$-scheme and reflection principles would be an interesting topic, but beyond the scope of this paper. 
	
	\begin{question}
		What is the exact consistency strength of $\mathsf{CZF+Sep}$ with a Reinhardt set?
	\end{question}
	Despite that, we will see in the next section that $\mathsf{CZF+Sep}$ with a Reinhardt set is still quite strong.
	
	\section{An analysis on the critical point}\label{Section:CriticalPoint}
	In the previous section, we saw that non-trivial elementary embeddings over $\mathsf{CZF+Sep}$ have a strong consistency strength. However, we observed nothing about how much large cardinal properties of the critical point of an $j$ are preserved by the double negation translation.
	In this section, we will extract the large cardinal properties of the critical point of an elementary embedding. Especially, we will focus on the critical point of a Reinhardt embedding.
	
	The main result of this section is as follows:
	\begin{theorem}\label{Theorem:ReflectionPropertyCriticalSet-weak}
		Let $j:V\prec M$ be a elementary embedding with an inaccessible critical point $K$. Then $V^\Omega$ thinks $\forall n,m\in\omega  (n<m\to \text{\normalfont $j^n(\tilde{K})$ is BCST-regular})^{j^m(\tilde{K})}$.
		
		Furthermore, if $j:V\prec V$, then $V^\Omega$ also thinks $j\restricts j^n(\tilde{K})\in j^{n+1}(\tilde{K})$.
	\end{theorem}
	
	The main strategy of this theorem is to internalize the proof of $V^\Omega\models \mathsf{CZF^-}$ into a BCST-regular set.
	We mostly follow the proof of \cite{Gambino2006}, but for the sake of verification, we will provide most of the detail of relevant lemmas and their proof. Throughout this section, $A\in B$ sets of the same power set of 1 (i.e., $\mathcal{P}(1)\cap A=\mathcal{P}(1)\cap B$) such that $A$ is BCST-regular and $B$ is a transitive model of $\mathsf{CZF^-}$, and $R\in B$ is a multi-valued function unless specified.
	
	\begin{lemma}\label{Lemma:TargetSubsetMVChoice}
		Let $a\in A$, $R:a\rrarrows A$, $Q\subseteq a\times A$, and $Q\in B$. Moreover, assume that
		\begin{enumerate}
			\item $\lag x,y\rag\in R\to y\subseteq Q_x=\{y\mid \lag x,y\rag\in Q\}$, and
			\item (Monotone Closeness) $\lag x,y\rag \in R$ and $y\subseteq z\subseteq Q_x$ implies $\langle x,z\rangle\in R$,
		\end{enumerate}
		then there is $f\in A\cap {^a}A$ such that $\lag x,f(x)\rag \in R$ for all $x\in a$.
	\end{lemma}
	\begin{proof}
		By \autoref{Lemma:SetMV}, there is $b\in A$ such that $b\subseteq R$ and $b:a\rrarrows A$. Take $f$ as $x\mapsto \bigcup b_x=\bigcup\{y\mid \lag x,y\rag \in b\}$.
		Then $f\in A$, since $A$ satisfies Union and the second-order Replacement.
		Moreover, by monotone-closedness of $R$, we have $\lag x,f(x)\rag\in R$ for all $x\in a$.
	\end{proof}
	
	\begin{lemma}\label{Lemma:TargetSubsetMVChoice-uniform}
		Let $P\subseteq A$, $P\in B$, $a\in A$ and $R:a\rrarrows A\cap\mathcal{P}(P)$. Furthermore, assume that $R$ is \emph{monotone closed}, that is, $y,z\in A$, $y\subseteq z\subseteq P$, and $\lag x,y\rag \in R$ implies $\lag x,z\rag \in R$. Then there is $b\in A$ such that $b\subseteq P$ and $\lag x,b\rag\in R$ for all $x\in a$.
	\end{lemma}
	\begin{proof}
		Applying \autoref{Lemma:TargetSubsetMVChoice} to $R$ provides a function $f\in A\cap {^a}A$ such that $\lag x,f(x)\rag \in R$ for all $x\in A$. Now take $b=\bigcup\{f(x)\mid x\in a\}\in A$, then we have $\lag x,b\rag\in R$ by monotone closeness of $R$.
	\end{proof}

	The following lemma has a critical role in the proof of our theorem. Moreover, this lemma requires $\mathcal{P}(1)\cap A=\mathcal{P}(1)\cap B$:
	\begin{lemma}\label{Lemma:Preserving-MainLemma}
		Let $a\in A^\Omega$, $R:a\rrarrows A$, $R\in B$ and $p\in A$ be such that $p\subseteq 1$ and $p=p^\dnot$.
		Define
		\begin{equation*}
			P = \{\lag x,y,t\rag\in \dom a\times A^\Omega\times 1 \mid t\in (p\land a(x)\land \llbr \op(x,y) \in R\rrbr^B)\}.
		\end{equation*}
		Furthermore, assume that we have $p\subseteq\llbr R:a\rrarrows \tilde{A}\rrbr^B$.
		Then there is $r\in A$ such that $r\subseteq P$ and $p\land a(x)\subseteq \{t\mid \exists y \in A^\Omega \lag x,y,t\rag \in r\}^\dnot$.
	\end{lemma}
	\begin{proof}
		Let $Q=\{\lag x,t\rag\mid \exists y\in A^\Omega (\lag x,y,t\rag\in P)\}$ and $Q_x=\{t\mid \lag x,t\rag \in Q\}\subseteq 1$.
		Then $Q_x\in \mathcal{P}(1)\cap B=\mathcal{P}(1)\cap A$. From $p\subseteq \llbr R:a\rrarrows \tilde{A}\rrbr^B$, we can deduce $p\land a(x)\subseteq Q_x^\dnot$.
		By applying \autoref{Lemma:TargetSubsetMVChoice} to the relation
		\begin{equation*}
			\{\lag x, v\rag \in\dom a\times (\mathcal{P}(1)\cap A)\mid p\land a(x)\subseteq v^\dnot \text{ and } v\subseteq Q_x\}
		\end{equation*}
		we have a function $f\in {^{\dom a}}A\cap A$ such that $p\land a(x)\subseteq f(x)^\dnot$ and $f(x)\subseteq Q_x$ for all $x\in\dom a$. (The condition $\mathcal{P}(1)\cap B=\mathcal{P}(1)\cap A$ is necessary to ensure the above relation is a multi-valued function of domain $\dom a$.)
		Now let
		\begin{equation*}
			q = \{\lag x,t\rag \mid x\in\dom a\text{ and } t\in f(x)\}.
		\end{equation*}
		Then $\forall \lag x,t\rag\in q\exists y\in A^\Omega (\lag x,y,t\rag \in P)$ holds. That is, $P:q\rrarrows A^\Omega$. By \autoref{Lemma:SetMV}, there is $r\in A$ such that $r\subseteq P$ and $r:q\rrarrows A^\Omega$.
		It is easy to see that $r$ satisfies the desired property.
	\end{proof}
	There is some technical note for the proof: there is no need that $P$, $Q$, and $Q_x$ are definable over $A$ in general. The reason is that we do not know either $R$ or $\llbr\cdot \rrbr^B$ is accessible from $A$. However, we do not need to care about it since we are relying on the second-order Strong Collection over $A$.
	
	\begin{theorem}\label{Theorem:BCSTset-preserving}
		Let $A$ be a BCST-regular set. Then $B^\Omega$ thinks $\tilde{A}$ is BCST-regular.
	\end{theorem}
	\begin{proof}
		It is easy to see that $B^\Omega$ thinks $\tilde{A}$ is transitive and closed under Pairing, Union, and Binary Intersection.
		Hence it remains to show that $B^\Omega$ thinks $\tilde{A}$ satisfies second-order Strong Collection, that is, 
		\begin{equation*}
			\llbr \forall a\in \tilde{A}\forall R [R:a\rrarrows \tilde{A}\to \exists b\in \tilde{A}(R:a\lrlrarrows b)] \rrbr^B = \top.
		\end{equation*}
		Take $a\in \dom\tilde{A}$, $R\in B^\Omega$ and $p\in A$ such that $p\subseteq 1$ and $p=p^\dnot$.
		We claim that if $p\subseteq \llbr R:a\rrarrows \tilde{A}\rrbr^B$, then there is $b\in \dom\tilde{A}$ such that $p\subseteq \llbr R:a\lrlrarrows b\rrbr^B$.
		By \autoref{Lemma:Preserving-MainLemma}, we have $r\in A$ such that $r\subseteq P$ and $p\land a(x)\subseteq \{t\mid \exists y \lag x,y,z\rag\in r\}^\dnot$.
		Define $b$ such that $\dom b=\{y\mid\exists x,t(\lag x,y,t\rag\in r)\}$ and 
		\begin{equation*}
			b(y) =\{0\mid \exists x \lag x,y,0\rag\in r\}^\dnot.
		\end{equation*}
		for $y\in\dom b$. (Note that $b(y)\in A$.) Then we have $p\subseteq \llbr R:a\lrlrarrows b\rrbr^B$.
	\end{proof}
	
	Hence we have
	\begin{corollary}\label{Corollary:reflectionK-jK-individual}
		$\llbr j^n(\tilde{K})\text{\text{ \normalfont is BCST-regular}}\rrbr^{j^m(K)}=\top$. Furthermore, $\llbr(j^n(\tilde{K})\text{\text{ \normalfont is BCST-regular}})^{j^m(\check{K})}\rrbr=\top$.
	\end{corollary}
	\begin{proof}
		The first statement follows from \autoref{Theorem:BCSTset-preserving} by taking $A=K$ and $B=j^m(K)$ for $n=1$, and applying $j$ $n-1$ times for general cases.
		For the last statement, we can apply the third clause of \autoref{Lemma:HeytingInterpretationBetweenModels}.
	\end{proof}
	However, it does not directly result in our desired theorem, since we do not know $\Omega$ forces $j^n(\tilde{K})$ is \emph{uniformly} BCST-regular. We need some work to see this.
	
	There are two possible meanings of $j^n(K)$: the first is applying $j$ $n$ times to $K$. Here $n$ must be a natural number over the metatheory, and this description lacks a way to describe the sequence $\langle j^n(K)\mid n\in\omega\rangle$. The second way to see $j^n(K)$ is to understand it as it is given by the following recursion:
	\begin{equation*}
		j^0(K)=K\text{ and } j^{n+1}(K)=j(j^n(K)).
	\end{equation*}
	Thus the formal statement of $\phi(j^n(K))$ is
	\begin{equation*}
		\exists f [\dom f=\omega\land f(0)=K\land \forall m\in\omega (f(m+1)=j(f(m)))]\land \phi(f(n)).
	\end{equation*}
	
	\begin{theorem}\label{Theorem:PropertyjnK-inVOmega}
		$V^\Omega$ thinks the following statement is valid:
		\begin{multline}
			\exists f [\text{$f$ is a function}\land \dom f=\omega\land f(0)=K\land \forall m\in\omega (f(m+1)=j(f(m)))]\\ \land
			\forall n,m\in\omega [n<m\to (f(n)\text{ \normalfont is BCST-regular})^{f(m)}].
		\end{multline}
	\end{theorem}
	\begin{proof}
		Let $\dom f=\{\op(\check{n},j^n(\tilde{K}))\mid n\in\omega\}$. We claim that $f$ witnesses our theorem. The first three conditions are easy to prove. 
		To see the last condition, let $n,m\in\omega$. We can show the following facts by induction on $m$:
		\begin{enumerate}
			\item If $0\in\llbr\check{n}<\check{m}\rrbr$ then $n<m$.
			\item $\llbr f(m)=j^m(\tilde{K})\rrbr=\top$.		
		\end{enumerate}
		By combining these facts with \autoref{Corollary:reflectionK-jK-individual}, we have $\llbr \check{n}<\check{m}\rrbr\subseteq \llbr(f(n)\text{ \normalfont is BCST-regular})^{f(m)}\rrbr$. Hence the result follows.
	\end{proof}

	Hence $\tilde{K}$ in $V^\Omega$ has the following reflection property: For every $n<m$, $j^m(\tilde{K})\models (j^n(\tilde{K})\models \mathsf{CZF^-_2})$.
	How much is this reflection principle strong? 
	To see this, assume that we started from $\mathsf{CZF+Sep}$ with an elementary embedding $j:V\prec V$ and produced the Heyting interpretation under $\Omega$. Then the Heyting interpretation interprets $\mathsf{ZF^-}$ with a cofinal elementary embedding $j:V\prec V$.
	By \autoref{Theorem:ReflectionPropertyCriticalSet-weak}, a critical point $\tilde{K}$ of $j$ satisfies $j(\tilde{K})\models (\tilde{K}\models \mathsf{CZF^-_2})$. Due to the help of the classical logic, we have $j(\tilde{K})\models (\tilde{K}\models \mathsf{ZF_2^-})$.
	For each $x\in \tilde{K}$, we have $j(\tilde{K})\models \exists X (x\in X\land X\models\mathsf{ZF_2^-})$.
	%
	%
	By the property of $j$, we have $\tilde{K}\models \forall x\exists X (x\in X\land X\models\mathsf{ZF^-_2})$. That is, $\tilde{K}$ satisfies $\mathsf{REA}$.
	Since $\mathsf{REA}$ implies the Axiom of Subset Collection, which is equivalent to Power Set in the classical context, we have $\tilde{K}\models \mathsf{ZF}$ and $j(\tilde{K})\models (\tilde{K}\models \mathsf{ZF_2})$!
	We may extend it further to stronger notions of large cardinal properties, like inaccessibility, Mahloness, or indescribability.
	However, the following example describes there is a limit of large cardinal property we can achieve from the reflection property of $\tilde{K}$:
	\begin{example}
		Work over $\mathsf{ZF}$ with a Reinhardt cardinal. Let $j:V\prec V$ be an elementary embedding and $\kappa=\crit j$. Take $K=L_\kappa$. Since $\kappa$ is a critical cardinal, it is strongly inaccessible by Proposition 3.3 of \cite{HayutKaragila2020}.
		Hence $L_{j(\kappa)}$ also thinks $\kappa$ is inaccessible. Especially, $L_{j(\kappa)}$ thinks $V_\kappa=L_\kappa$ is a model of $\mathsf{ZFC_2}$.
	\end{example}
	Since $L$ is incompatible with large cardinals stronger than the existence of $0^\sharp$, the above example shows the previous argument with $K$ does not yield large cardinal properties stronger than the existence of $0^\sharp$.
	However, it does not mean there is no room for stronger properties of $K$ if $j$ has a stroger property. The following result shows we can extract more large cardinal properties from $K$ if $j:V\prec V$, by bringing the elementary embedding $j$ over $\mathsf{CZF}$ to an elementary embedding of $j^\omega(\tilde{K})$ in the Heyting interpretation:
	\begin{lemma}[$\mathsf{CZF^-}$]
		Let $j:V\prec V$. Then
		$V^\Omega$ thinks $j\restricts j^n(\tilde{K})\in j^{n+1}(\tilde{K})$ for all $n\in\omega$.
	\end{lemma}
	\begin{proof}
		Observe that $j^{n+1}(K)$ is regular and $j\restricts j^n(K):j^n(K)\to j^{n+1}(K)$ is a multi-valued function. By \autoref{Lemma:SetMV}, there is $b\in j^n(K)$ such that $b\subseteq j\restricts j^n(K)$ and $b:j^n(K)\rrarrows j^{n+1}(K)$. Since $j\restricts j^n(K)$ is a function, $b$ is also a function of domain $j^n(K)$. Hence $j\restricts j^n(K)=b\in j^{n+1}(K)$.
		
		From the previous argument, we also have $j\restricts j^n(\tilde{K})=(j\restricts j^n(K))\restricts j^n(\tilde{K})\in j^{n+1}(K)$. Now let $c$ be a name such that
		\begin{equation*}
			\dom c=\{\op(x,y)\mid \langle x,y\rangle \in j\restricts j^n(\tilde{K})\}
			\quad\text{and}\quad
		\end{equation*}
		and $c(x)=\top$ for all $x\in \dom c$. By definition of $c$, $c\in j^{n+1}(\tilde{K})$.
		Moreover, it is easy to see that $V^\Omega$ thinks $c$ is a function of the domain $j^n(\tilde{K})$, and direct calculation shows $\llbr\forall x\forall y \op(x,y)\in c\to j(x)=y\rrbr=\top$. Hence $\llbr c = j\restricts j^n(\tilde{K})\rrbr=\top$.
	\end{proof}
	Note that this lemma does not work for general $j:V\prec M$, since we do not know $j^n(K)$ is  regular in $V$.
	
	\begin{theorem}[$\mathsf{CZF+Sep}$]\label{Theorem:AnalysisEmbedding-main}
		$V^\Omega$ thinks $j^\omega(\tilde{K}):=\bigcup_{n\in\omega} j^n(\tilde{K})$ satisfies $\mathsf{ZF}$ with a cofinal elementary embedding $j$ from itself to itself. Moreover, $j^\omega(\tilde{K})$ satisfies $\Delta_0$-Separation with $j$ be allowed to appear. 
	\end{theorem}
	\begin{proof}
		Let $\tilde{K}=j(\tilde{K})$ and $\Lambda:=j^\omega(\tilde{K})$. We will rely on a completely internal argument to $V^\Omega$, which is a model of $\mathsf{ZF^-}$.
		
		Assume that $\Lambda\models \phi(\vec{a})$, where $\phi(\vec{x})$ is a $\Delta_0$-formula with all free variables displayed in the language $\in$ (i.e., without $j$.) Since $j(\Lambda)=\Lambda$, we have $\Lambda\models \phi(j(\vec{a}))$. Hence $j:V\prec V$.
		
		It remains to show that $\Lambda$ satisfies $\Delta_0$-separation for formulas with $j$ be allowed to appear.
		Let $a\in\Lambda$, then there is $n$ such that $a\in j^n(\tilde{K})$. For a bounded formula $\phi$ with parameters in $j^n(\tilde{K})$, let $\phi'$ be the formula obtained by every occurrence of $j$ to $j\restricts j^n(\tilde{K})$. Then $\{x\in a\mid \phi(x)\}=\{x\in a\mid \phi'(x)\}\in j^{n+1}(\tilde{K})$. Thus $\Lambda$ satisfies $\Delta_0$-separation for formulas with $j$.
	\end{proof}
	
	How much is the resulting theory strong? We can see that $\Lambda$ is a model of $\mathsf{ZF}$ with the Wholeness axiom for $\Delta_0$-formulas $\mathsf{WA_0}$ proposed by Hamkins \cite{Hamkins2001WA}. The author does not know the exact consistency strength of $\mathsf{ZF+WA_0}$ in $\mathsf{ZFC}$-context. However, we can still find a lower bound of it: we can see that $\Lambda$ also thinks $\kappa=\rank K$ is a critical point of $j$, and the critical sequence defined by $\kappa_0=\kappa$ and $\kappa_{n+1}=j(\kappa_n)$ is cofinal over $\Ord$. (Note that the cofinal sequence is still definable, although it may not a set. See Proposition 3.2 of \cite{Corazza2000} for details.) From this, we have
	\begin{lemma}[$\mathsf{ZF+WA_0}$]
		If the critical sequence is cofinal, then $\kappa_0$ is extendible.
	\end{lemma} 
	\begin{proof}
		Let $\eta$ be an ordinal. Take $n$ such that $\eta<j^n(\kappa)$, then $j^n:V_{\kappa+\eta}\prec V_{j^n(\kappa+\eta)}$ and $\crit j^n=\kappa$. Hence $\kappa$ satisfies $\eta$-extendibility.
	\end{proof}
	By an easy reflection argument, we can see also see that $\mathsf{ZF+WA_0}$ with the cofinal critical sequence proves not only there is an extendible cardinal, but also the consistency of $\mathsf{ZF}$ with a proper class of extendible cardinals, an extendible limit of extendible cardinals, and many more.
	Since extendible cardinals are preserved by Woodin's forcing (Theorem 226 of \cite{Woodin2010}), we have a lower bound of the consistency strength of $\mathsf{ZF+WA_0}$, e.g., $\mathsf{ZFC}$ with there is a proper class of extendible cardinals.
	
	\section{Concluding Questions}\label{Section:RemarkQuestions}
	We may wonder how to find the upper bound of the consistency strength of $\mathsf{CZF+Sep}$ with a non-trivial elementary embedding, in terms of extensions of $\mathsf{ZFC}$ or $\mathsf{ZFC^-}$.
	Most construction of interpretations of $\mathsf{CZF}$ from classical theories rely on realizability, and employ type-theoretic interpretations (like \cite{Rathjen2005Brouwer} or \cite{GrifforRathjen1994}) or set-as-tree interpretation (like \cite{Lubarsky2006SOA} or functional realizability model over $\mathsf{ZFC^-}$ given by Swan \cite{Swan2014}.\footnote{Swan worked over $\mathsf{ZFC}$, but his proof for soundness of functional realizability still works over $\mathsf{ZFC^-}$. Also, note that his functional realizability is a part of his two-stage Kripke model.})
	These interpretations usually satisfy the \emph{Axiom of Subcountability}, which states every set is an image of a subset of $\omega$. However, Ziegler \cite{ZieglerPhD} proved that the existence of non-trivial elementary embedding $j:V\prec M$ contradicts with the Axiom of Subcountability. The author thinks that the Axiom of Subcountability comes from that we are using countable pca to construct interpretations, but delimiting the size of pca does not guarantee that we can reach the upper bound of $\mathsf{CZF}$ with a non-trivial elementary embedding.
	We do not know about the type-theoretic analogue of $\mathsf{CZF}$ with an elementary embedding, so we do not know we can employ type-theoretic interpretations. 
	Functional realizability or set-as-tree interpretations also have a problem. We need a pca of size greater than not only the critical point but also its successive application to the elementary embedding to ensure the critical set exists under the interpretation since the size of the pca delimits the size of every set under the interpretation. However, nothing much is known about realizability under large pcas. Especially, these large pcas are not fixed under the elementary embedding, which makes them hard to handle. One possible way to control the large pcas under the elementary embedding $j$ is to add $j$ into the pca. However, finding the realizer of $\phi^M(j(x))\to\phi(x)$ seems not easy.
	
	\begin{question}
		What is the upper bound of the consistency strength of $\mathsf{CZF+Sep}$ with a critical set or a Reinhardt set? Can we find the bound in terms of large cardinals compatible with the Axiom of Choice?
	\end{question}
	
	Our method is also restricted to the analysis on $\mathsf{CZF+Sep}$, mainly because the forcing over the double negation topology only provides $\mathsf{\Delta_0\mhyphen LEM}$. Full Separation is necessary to turn it to the full excluded middle. We may ask we can analyze the strength of $\mathsf{CZF}$ with large large set axioms without any help of full Separation. 
	\begin{question}
		Is there any non-trivial result for the consistency strength of $\mathsf{CZF}$ with a critical set or a Reinhardt set?
	\end{question}
	
	Improving the lower bound of the consistency strength of $\mathsf{CZF+Sep}$ with a Reinhardt set is also an issue. We have not made use of the full strength of the resulting theory stated in \autoref{Proposition:MainResult} to get the lower bound of the consistency strength.
	For example, we never used the cofinality of $j$. Even worse, we did not extract the full consistency strength of $\mathsf{ZF+WA_0}$: extendibility (or a proper class of extendibles) is far below $\mathsf{WA_0}$.
	One may try to force the Axiom of Choice over $\mathsf{ZF+WA_0}$, by using Woodin's forcing (see Theorem 226 of \cite{Woodin2010}) and a method given by Hamkins \cite{Hamkins2001WA}. However, the quotient of the limit stages of Woodin's forcing is not sufficiently $\gamma$-closed, so Hamkins' argument is not applied well.
	\begin{question}
		Can we obtain a better lower bound of the consistency strength of $\mathsf{ZF^-}$ with a cofinal embedding $j:V\prec V$, with a critical point $K\models \mathsf{ZF}$ such that $j^{\omega}(K)\models (K\models \mathsf{ZF_2^-})$?
		Especially, are $\mathsf{ZF+WA_0}$ and $\mathsf{ZFC+WA_0}$ equiconsistent?
	\end{question}
	
	\section*{Acknowledgements}
	The author wants to thank Richard Matthews for pointing out an error and providing comments on earlier versions of this manuscript.

	\overfullrule=0pt
	\printbibliography
	
\end{document}